\documentclass[12pt]{article}
\usepackage{amsfonts, amsmath}
\usepackage{amsfonts}
\usepackage{amssymb}
\usepackage{amsfonts}

\newtheorem{theorem}{Theorem}[section]

\newtheorem{condition}[theorem]{Condition}

\newtheorem{definition}[theorem]{Definition}
\newtheorem{example}[theorem]{Example}

\newtheorem{proposition}[theorem]{Proposition}
\newtheorem{remark}[theorem]{Remark}

\newenvironment{proof}[1][Proof]{\noindent\textbf{#1.} }{\ \rule{0.5em}{0.5em}}
\begin{document}

\title{Stratified Subcartesian Spaces}
\author{Tsasa Lusala and J\c{e}drzej \'{S}niatycki\\
\\
Department of Mathematics and Statistics\\
University of Calgary\\
2500 University Drive N.W.\\
Calgary, Alberta T2N1N4, Canada\\
Email: $\{$tsasa, sniat$\}$@math.ucalgary.ca}

\date{}
\maketitle

\begin{abstract}
We show that, if the family $\mathcal{O}$ of orbits of all vector fields on
a subcartesian space $P$ is locally finite and each orbit in $\mathcal{O}$
is locally closed, then $\mathcal{O}$ defines a smooth Whitney A
stratification of $P$. We also show that the stratification by orbit type of
the space $M/G$ of orbits of a proper action of a Lie group $G$ on a smooth
manifold $M$ is given by orbits of the family of all vector fields on $M/G$.
\end{abstract}

\noindent
{\bf 2000 Mathematics Subject Classification: 58A40, 57N80.}

\medskip

\noindent
{\bf Keywords and Phrases}: Subcartesian spaces, orbits of vector fields, Stratifications, Whitney Conditions.

\section{Introduction}

Stratification theory is based on the natural idea of dividing a singular
space into manifolds. It deals with study of topological spaces endowed with
a partition by smooth manifolds satisfying specific conditions. Many of
singular spaces appearing in analysis have the structure of \ stratified
spaces satisfying Whitney's condition B \cite{whitney}, and the theory of
stratified spaces is an important tool with a broad range of applications,
see \cite{goresky-macpherson} and references quoted there.

Sikorski's theory of differential spaces is a tool in the study of the
differential geometry of a large class of singular spaces \cite{sikorski}. A
differential space $P$ is said to be subcartesian if every point $p\in P$
has a neighbourhood diffeomorphic to a subset of a Euclidean space \cite%
{aronszajn}. In particular, an arbitrary subset $P$ of $\mathbb{R}^{n}$,
with the ring of smooth functions generated by restrictions to $P$ of smooth
functions on $\mathbb{R}^{n}$, is subcartesian.

Every subcartesian space has a canonical partition by smooth manifolds given
by orbits of the family of all vector fields on the space \cite{sniatycki}.
The aim of this paper is to discuss stratifications of subcartesian spaces
and compare them with partitions by orbits of the family of all vector
fields. We show that the partition of a subcartesian space $P$ by the family 
$\mathcal{O}$\ of orbits all vector fields satisfies the frontier condition
and Whitney's condition A. From this we conclude that, if the family $%
\mathcal{O}$ is locally finite and each orbit in $\mathcal{O}$ is locally
closed, then $\mathcal{O}$ defines a smooth Whitney A stratification of $P$.
A locally finite family $\mathcal{O}$ of locally closed orbits of all vector
fields need not satisfy Whitney's condition B. However, some smooth Whitney
B stratifications are given by orbits of all vector fields. We show that the
stratification by orbit type of the space of orbits $M/G$ of a proper action
of a Lie group $G$ on a smooth manifold $M$ is given by orbits of the family
of all vector fields on $M/G$.

\section{Decomposed spaces}

A decomposition of a topological space $P$ is a partition of $P$ by a
locally finite family $\mathcal{D}$ of smooth manifolds $M\subset P$, such
that each manifold $M\in \mathcal{D}$ with its manifold topology is a
locally closed topological subspace of $P$, satisfying the following
condition:

\begin{condition}[Frontier Condition]
For $M,M^{\prime }\in \mathcal{D}$, if $M^{\prime }\cap \bar{M}\neq
\emptyset $, then either $M^{\prime }=M$ or $M^{\prime }\subset \bar{M}%
\backslash M$.
\end{condition}

\noindent The pair $(P,\mathcal{D})$ is called a decomposed space. Local
finiteness of $\mathcal{D}$ means that, for each point $p\in P$, there
exists a neighbourhood $U$ of $p$ in $P$ intersecting only a finite number
of manifolds $M\in \mathcal{D}$. A subset $M$ of a topological space $P$ is
locally closed if, for each $x\in M$ there exists a neighbourhood $U$ of $x$
in $P$ such that $M\cap U$ is closed in $U.$ If $P$ is a manifold, an
injectively immersed submanifold $M$ of $P$ is embedded if and only if $M$
is locally closed in $P.$

Decomposed spaces form a category with morphisms $\varphi :(P_{1},\mathcal{D}%
_{1})\rightarrow (P_{2},\mathcal{D}_{2})$ given by continuous map $\varphi
:P_{1}\rightarrow P_{2}$ such that, for each $M_{1}\in \mathcal{D}_{1},$
there exists $M_{2}\in \mathcal{D}_{2}$ such that $\varphi (M_{1})\subset
M_{2}$ and the restriction of $\varphi $ to $M_{1}$ is a smooth map from $%
M_{1}$ to $M_{2}$.

For a decomposed space $(P,\mathcal{D})$, let $Q$ be a topological subspace
of $P$, and $\mathcal{D}_{Q}=\{M\cap Q\mid M\in \mathcal{D}\}$. Suppose
that, for each $M\in \mathcal{D}$, $M\cap Q$ is a submanifold of $M$ locally
closed in $Q$, and the family $\mathcal{D}_{Q}$ is locally finite. Then, $%
\mathcal{D}_{Q}$ satisfies the Frontier Condition because, if $M$ and $%
M^{\prime }$ in $\mathcal{D}$ are such that $(M^{\prime }\cap Q)\cap 
\overline{(M\cap Q)}\neq \emptyset $, then $M^{\prime }\cap \bar{M}\neq
\emptyset $, and either $M^{\prime }=M$ or $M^{\prime }\subset \bar{M}-M$,
so that, either $(M^{\prime }\cap Q)=(M\cap Q)$ or $(M^{\prime }\cap
Q)\subset (M\cap Q)$. Therefore, $(Q,\mathcal{D}_{Q})$ is a decomposed
space. In particular, if $Q$ is an open subset of $P$, then $(U,\mathcal{D}%
_{U})$ is a decomposed space.

Suppose $(P,\mathcal{D})$ is a decomposed space, $Q$ is a smooth manifold,
and $\mathcal{D}_{P\times Q}=\{M\times Q\mid M\in \mathcal{D}\}$. Then $%
(P\times Q,\mathcal{D}_{P\times Q})$ is also a decomposed space, and the
projection map $P\times Q\rightarrow P$ gives a morphism from $(P\times Q,%
\mathcal{D}_{P\times Q})$ to $(P,\mathcal{D})$. A decomposed space $(P,%
\mathcal{D})$ is locally trivial if, for every point $M\in \mathcal{D}$ and
each $x\in M$, there exists an open neighbourhood $U$ of $x$ in $P$, a
decomposed space $(P^{\prime },\mathcal{D}^{\prime })$ with a distinguished
point $y\in P^{\prime }$ such that the singleton $\{y\}\in \mathcal{D}%
^{\prime }$, and an isomorphism $\varphi :(U,\mathcal{D}_{U})\rightarrow
(P^{\prime }\times (U\cap M),\mathcal{D}_{P^{\prime }\times (U\cap
M)}^{\prime }),$ such that $\varphi (x)=y$.

Decompositions of a topological space $P$ can be partially ordered by
inclusion. If $\mathcal{D}_{1}$ and $\mathcal{D}_{2}$ are two decompositions
of $P$, we say that $\mathcal{D}_{1}$ is a \emph{refinement} of $\mathcal{D}%
_{2}$ and write $\mathcal{D}_{1}\geq \mathcal{D}_{2}$, if, for every $%
M_{1}\in \mathcal{D}_{1}$ there exists $M_{2}\in \mathcal{D}_{2}$ such that $%
M_{1}\subseteq M_{2}$. We say that $\mathcal{D}$ is a minimal (coarsest)
decomposition of $P$ if it is not a refinement of a different decomposition
of $P$. Note that if $P$ is a manifold, then the minimal decomposition of $M$
consists of a single manifold $M=P$. Similarly, we say that $\mathcal{D}$ is
a maximal (finest) decomposition of $P$ if $\mathcal{D}^{\prime }\geq 
\mathcal{D}$ implies $\mathcal{D}^{\prime }=\mathcal{D}$.

\section{Stratified spaces}

Let $A$ and $B$ be subsets of a topological space $P$. If $x\in A\cap B$, we
say that $A$ and $B$ are equivalent at $x$ if there exists a neighbourhood $%
U $ of $x$ in $P$ such that $A\cap U=B\cap U$. The equivalence class at $x$
of a subset $A$ of $P$ containing $x$ is called the germ of $A$ at $x$ and
denoted $[A]_{x}.$

A stratification of a topological space $P$ is a map $\mathcal{S}$ which
associates to each $x\in P$ a germ $\mathcal{S}_{x}$ of a manifold embedded
in $P$ such that the following condition is satisfied:

\begin{condition}[Stratification Condition]
For every $z\in P$ there exists a neighbourhood $U$ of $z$ and a
decomposition $\mathcal{D}$ of $U$ such that for all $y\in U$ the germ $%
\mathcal{S}_{y}$ coincides with the germ of the manifold $M\in \mathcal{D}$
that contains $y$.
\end{condition}

\noindent Every decomposition $\mathcal{D}$ of $P$ defines a stratification $%
\mathcal{S}$ of $P$ that associates to every $x\in P$ the germ $\mathcal{S}%
_{x}$ at $x$ of the manifold $M\in \mathcal{D}$ that contains $x$.

\begin{definition}
Two decompositions $\mathcal{D}_{1}$ and $\mathcal{D}_{2}$ of $P$ are
equivalent if they define the same stratification $\mathcal{S}$ of $P$.
\end{definition}

Let $\mathcal{S}$ be a stratification of $P$. There is a unique
decomposition $\mathcal{D}_{0}$ of $P$ by connected manifolds which defines $%
\mathcal{S}$. It is the finest element of the class of decomposition of $P\ $%
corresponding to $\mathcal{S}$. From the point of view of this paper it is
convenient to identify $\mathcal{S}$ with $\mathcal{D}_{0}$.

\section{Differential and subcartesian spaces}

A differential structure on a topological space $P$ is a family $C^{\infty
}(P)$ of functions on $P$ satisfying the following conditions.

\begin{condition}[Differential Structure]
\begin{enumerate}
\item The family of sets $\{f^{-1}((a,b))\mid f\in C^{\infty }(P)$, and $%
a,b\in \mathbb{R}\}$ is a sub-basis for the topology of $P$.

\item For every $k\in \mathbb{N}$, every $f_{1},...,f_{k}\in C^{\infty }(P)$
and $F\in C^{\infty }(\mathbb{R}^{k})$, the composition $F(f_{1},...,f_{k})$
is in $C^{\infty }(P)$.

\item If a function $f$ on $P$ is such that, for every $x\in P$, there
exists an open neighbourhood $U_{x}$ of $x$ in $Q$, and a function $f_{x}\in
C^{\infty }(P)$ satisfying $f_{\mid U_{x}}=f_{x\mid U_{x}},$ then $f$ is in $%
C^{\infty }(P)$.
\end{enumerate}
\end{condition}

\noindent A topological space endowed with a subring of continuous functions
satisfying the above conditions is called a differential space. 

A homeomorphism $\varphi :P\rightarrow Q$ of differential spaces is smooth
if its pull-back $\varphi ^{\ast }$ maps $C^{\infty }(Q)$ to $C^{\infty }(P)$%
. It is a diffeomorphism if it is invertible and $\varphi ^{-1}:Q\rightarrow
P$ is smooth. A subcartesian space is a Hausdorff differential space $P$
such that, each point $x\in P$ has a neigbourhood which is diffeomorphic to
a subset of a Cartesian space $\mathbb{R}^{N}$.

We can adapt notions of decomposition and stratification of a topological
space to a differential space by requiring smoothness of all maps involved.
Thus, a \emph{smooth decomposition} of a differential space $P$ is a
decomposition $\mathcal{D}$ of $P$ as a topological space such that, for
each $M\in \mathcal{D}$, the inclusion map $M\hookrightarrow P$ is smooth.
Similarly, a \emph{smooth stratification} of a differential space $P$ is a
smooth decomposition of $P$ by connected manifolds.

For each point $p$ of a differential space $P$, a derivation of $C^{\infty
}(P)$ at $p$ is a linear map $u:C^{\infty }(P)\rightarrow \mathbb{R}$
satisfying Leibniz' rule 
\[
u(f_{1}f_{2})=f_{1}(p)u(f_{2})+u(f_{1})f_{2}(p)\text{ for all }%
f_{1},f_{2}\in C^{\infty }(P).
\]%
The set of all derivations of $C^{\infty }(P)$ at $p$ is called the \emph{%
tangent space} to $P$ at $p.$ It is denoted $T_{p}P$. The \emph{tangent cone}
at $p\in P$ is the subset $T_{p}^{C}P$ of $T_{p}P$ consisting of derivations
at $p$ which are given by differentiation along smooth curves in $P$ passing
through $p$. In other words, $u\in T_{p}P$ is in $T_{p}^{C}P$ if there
exists a smooth curve $c:[0,1]\rightarrow P$ such that $u(f)=\frac{d}{dt}%
f(c(t))_{\mid t=0}$. Reparametrization of curves gives rise to the cone
structure in $T_{p}^{C}(P).$

A (global) \emph{derivation} of $C^{\infty }(P)$ is a linear map $%
X:C^{\infty }(P)\rightarrow C^{\infty }(P)$ satisfying Leibniz' rule%
\[
X(f_{1}f_{2})=f_{1}X(f_{2})+X(f_{1})f_{2}\text{ for all }f_{1},f_{2}\in
C^{\infty }(P).
\]%
Let $I$ be an interval in $\mathbb{R}$ with a non-empty interior. A smooth
map $c:I\rightarrow P$ is an integral curve of a derivation $X$ if, for
every $t\in I$, and $f\in C^{\infty }(P),$ 
\[
\frac{d}{dt}f(c(t))=X(f)(c(t)).
\]%
The notion of an integral curve can be extended to the case when $I$
consists of a single point, i.e. $I=[a,a]$ for $a\in \mathbb{R}.$ In this
case the left hand side of the above equation is not defined. We consider a
map $c:I=[a,a]\rightarrow P:a\mapsto c(a)$ to be an integral curve of every
derivation of $C^{\infty }(P)$. With this definition, for every derivation $X
$ of the differential structure $C^{\infty }(P)$ of a subcartesian space $P$
and every $x\in P$, there exists a unique maximal integral curve of $X$
passing through $x$ \cite{sniatycki}.

\begin{definition}
A derivation $X$ of $C^{\infty }(P)$ is a vector field on a subcartesian
space $P$ if translations along integral curves of $X$ give rise to a local
one-parameter group of local diffeomorphisms of $P$.
\end{definition}

Let $\mathcal{X}$ be the family of all vector fields on $P$. For each $X\in 
\mathcal{X}$, we denote by $\exp tX$ the local one-parameter group of local
diffeomorphisms generated by translations along integral curves of $X$. The
orbit of $\mathcal{X}$ through a point $x\in P$ is 
\[
O_{x}=\{(\exp t_{1}X_{1}\circ ...\circ \exp t_{n}X_{n})(x)\mid n\in \mathbb{N%
},\text{ }(t_{1},...,t_{n})\in \mathbb{R}^{n},\text{ }X_{1},...,X_{n}\in 
\mathcal{X}\}. 
\]%
For each $x\in P$, the orbit $O_{x}$ of the family $X$ of all vector fields
on $P$ is a manifold and the inclusion map $O_{x}\hookrightarrow P$ is
smooth \cite{sniatycki}. The collection $\mathcal{O}$ of all orbits of $%
\mathcal{X}$ is a partition of $P$ by smoothly included manifolds.

\begin{theorem}
\label{theorem1} The partition $\mathcal{O}$ of a subcartesian space $P$ by
orbits of the family of all vector fields on $P$ is a smooth stratification
of $P$ if $\mathcal{O}$ is a locally finite, and each orbit $O\in \mathcal{O}
$ is locally closed.
\end{theorem}

\begin{proof}
By definition, orbits of the family of all vector fields are connected.
Moreover, for each orbit $O\in \mathcal{O},$ the inclusion map $%
O\hookrightarrow P$ is smooth. Hence, it suffices to show that the family $%
\mathcal{O}$ satisfies the Frontier Condition. Suppose $x\in O^{\prime }\cap 
\bar{O}$. with $O^{\prime }\neq O$. We first show that $O^{\prime }\subset 
\bar{O}$. Note that the orbit $O$ is invariant under the family of
one-parameter local groups of local diffeomorphisms of $P$ generated by
vector fields. Since, $x\in \bar{O}$, it follows that, for every vector
field $X$ on $S$, $\exp (tX)(x)$ is in $\bar{O}$ if it is defined. But, $%
O^{\prime }$ is the orbit of $\mathcal{X}$ through $x$. Hence, $O^{\prime
}\subset \bar{O}$.
\end{proof}

\smallskip 

A smoothly decomposed diffferential space $(P,\mathcal{D})$ is smoothly
locally trivial if, for every point $M\in \mathcal{D}$ and each $x\in M$,
there exists an open neighbourhood $U$ of $x$ in $P$, a smoothly decomposed
differential space $(P^{\prime },\mathcal{D}^{\prime })$ with a
distinguished point $y\in P^{\prime }$ such that the singleton $\{y\}\in 
\mathcal{D}^{\prime }$, and an isomorphism $\varphi :(U,\mathcal{D}%
_{U})\rightarrow (P^{\prime }\times (U\cap M),\mathcal{D}_{P^{\prime }\times
(U\cap M)}^{\prime }),$ such that $\varphi (x)=y$. It should be noted that a
smoothly decomposed differential space $(P,\mathcal{D})$ may be locally
trivial as a (topological) decomposed space but not smoothly locally
trivial. The following example, taken from Mather \cite{mather} illustrates
this situation.

\begin{example}
Consider $F(x,y,z)=xy(x+y)(x+\alpha (x)y)$ for a smooth one-to-one function $%
\alpha (x)$ with values different from $0$ and $1$. The zero level $S$ of $F$%
, given \ by 
\[
S=\{(x,y.z)\in R^{3}\mid xy(x+y)(y+\alpha (z)x)=0\}, 
\]%
is the union of four surfaces intersecting along the $z$-axis. It has eight
2-dimensional strata: $\pm x>0$, $\pm y>0$, $\pm (x+y)>0$ and $\pm (y+\alpha
(z))>0$, and a $1$-dimensional stratum consisting of the $z$-axis.

For each $z_{0},$ the tangent cone to $S$ at $(0,0,z_{0})$ is the union of
four planes $x=0$, $y=0$, $x+y=0$ and $x+\alpha (z_{0})y=0$ intersecting
along the $z$-axis. Projections of these planes to the $(x,y)$-plane are
four lines intersecting at the origin. If values of $\alpha $ are different
from $0$ and $-1$, then all four lines are distinct and their cross-ratio is 
$\gamma (z_{0})=1+\alpha (z_{0})$. By assumption, the function $\alpha (z)$
is one-to-one. Since the cross-ratio is an invariant of linear
transformations preserving the origin, all diffeomorphisms of $S$ preserve
points on the $z$-axis. This implies that the stratification of $S$
described above is not locally trivial.

The argument above implies also that the $z$-axis is not an orbit of the
family of all vector fields on $S$. The partition of $S$ by the family of
orbits of all vector fields on $S$ consists of 2-dimensional orbits, which
coincide with two dimensional strata, and 0-dimensional orbits $%
\{(0,0,z_{0})\}$ for each $z_{0}\in \mathbb{R}$.
\end{example}

Let $\mathcal{D}$ be a smooth decomposition of a differential space $P$. We
say that $\mathcal{D}$ admits \emph{local extension of vector fields} if,
for each $M\in \mathcal{D}$ each vector field $X_{M}$ on $M$ and each point $%
x\in M$, there exists a neighbourhood $V$ of $x$ in $M$, and a vector field $%
X$ on $P$ such that $X_{\mid V}=X_{M\mid V}$. In other words, the vector
field $X$ is an extension to $P$ of the restriction of $X_{M}$ to $V$.

\begin{theorem}
\label{theorem2}Every smoothly locally trivial decomposition of a
subcartesian space $P$ admits local extension of vector fields.
\end{theorem}

\begin{proof}
Let $X_{M}$ be a vector field on $M\in D$. Given $x_{0}\in M$, let $U$ be a
neighbourhood of $x_{0}$ in $M$ admitting an isomorphism $\varphi
:U\rightarrow (M\cap U)\times P^{\prime }$, for some smoothly decomposed
differential space $(P^{\prime },\mathcal{D}^{\prime })$ such that $%
\{\varphi (x_{0})\}\in \mathcal{D}^{\prime }$. Let $\exp (tX_{M})$ be the
local one-parameter group of local diffeomorphisms of $M$ generated by $X_{M}
$ and $X_{(M\cap U)\times P^{\prime }}$ be a derivation of $C^{\infty
}((M\cap U)\times P^{\prime })$ defined by 
\[
(X_{(M\cap U)\times P^{\prime }}h)(x,y)=\frac{d}{dt}h(\exp
(tX_{M})(x),y)_{\mid t=0}.
\]%
for every $h\in C^{\infty }((M\cap U)\times P^{\prime })$ and each $(x,y)\in
(M\cap U)\times P^{\prime }$. Since $X_{(M\cap U)\times P^{\prime }}$ is
defined in terms of a local one-parameter group $(x,y)\mapsto (\exp
(tX_{M})(x),y)$ of diffeomorphisms, it is a vector field on $(M\cap U)\times
P^{\prime }$.

We can use the inverse of the diffeomorphism $\varphi :U\rightarrow (M\cap
U)\times P^{\prime }$ to push-forward $X_{(M\cap U)\times P^{\prime }}$ to a
vector field $X_{U}=\varphi _{\ast }^{-1}X_{(M\cap U)\times P^{\prime }}$ on 
$U$. Choose a function $f_{0}\in C^{\infty }(P)$ with support in $U$ and
such that $f(x)=1$ for $x$ in some neighbourhood $U_{0}$ of $x_{0}$
contained in $U$. Let $X$ be a derivation $X$ of $C^{\infty }(P)$ extending $%
f_{0}X_{U}$ by zero outside $U$. In other words, for every $f\in C^{\infty
}(P)$, if $x\in U$, then $(Xf)(x)=f_{0}(x)(X_{U}f)(x),$ and if $x\notin
U_{0},$ then $(Xf)(x)=0$. Clearly, $X$ is a vector field on $P$ extending
the restriction of $X_{M}$ to $M\cap U_{0}$.
\end{proof}

\begin{theorem}
\label{theorem3}Let $\mathcal{D}$ be a decomposition of a subcartesian space 
$P$ admitting local extensions of vector fields, then the partition $%
\mathcal{O}$ of $P$ by orbits of the family of all vector fields on $P$ is a
stratification of $P$. If all manifolds in $\mathcal{D}$ are connected, then 
$\mathcal{D}$ is a refinement of $\mathcal{O}$. Moreover, if $\mathcal{D}$
is minimal in the class of decompositions by connected manifolds then $%
\mathcal{D}=\mathcal{O}$.
\end{theorem}

\begin{proof}
Let $\mathcal{D}$ be a decomposition of $P$ admitting local extension of
vector fields. Since every vector field $X_{M}$ on a manifold $M\in \mathcal{%
D}$ extends locally to a vector field on $P$, it follows that $M$ is
contained in an orbit $O\in \mathcal{O}$.

Every orbit $O\in \mathcal{O}$ is a union of manifolds in the decomposition $%
\mathcal{D}$. Since $\mathcal{D}$ is locally finite, it follows that, for
each $x\in P$, there exists a neigbourhood $U$ of $x$ in $P$ which
intersects only a finite number of manifolds in $\mathcal{D}$. Hence, $U$
intersects only a finite number of orbits in $\mathcal{O}$.

Since manifolds in $\mathcal{D}$ are locally closed, for each $M\in \mathcal{%
D},$ and each $x\in M$, there exists a neighbourhood $U$ of $x$ in $P$ such
that $M\cap U$ is closed in $U$. Without loss of generality, we may assume
that there is only a finite number of manifolds $M_{1}=M$, $M_{2}$, ..., $%
M_{k}$ in $\mathcal{D}$ such that $M_{i}\cap U\neq \emptyset $ for $%
i=1,...,k $. Since manifolds in $\mathcal{D}$ form a partition of $P$, it
follows that 
\[
U=\bigcup_{i=1}^{k}M_{i}\cap U. 
\]%
We may also assume that each $M_{i}\cap U$ is closed in $U$.

Let $O$ be the orbit in $\mathcal{O}$ that contains $M=M_{1}$. We can
relabel the manifolds $M_{1},...,M_{k}$ so that 
\[
O\cap U=O\cap \bigcup_{i=1}^{k}M_{i}\cap
U=\bigcup_{i=1}^{k}O\cap M_{i}\cap
U=\bigcup_{i=1}^{l}M_{i}\cap U 
\]%
for some $l\leq k$. Since $M_{i}\cap U$ is closed in $U$ for each $i=1,...,l$%
, it follows that $O\cap U$ is also closed in $U$. Hence, orbits $O\in 
\mathcal{O}$ are locally closed.

Taking into account Theorem \ref{theorem1} we see that $\mathcal{O}$ is a
stratification of $P$. If all manifolds in $\mathcal{D}$ are connected then
each $M\in \mathcal{D}$ is contained in an orbit in $\mathcal{O}$ and $%
\mathcal{D}$ is a refinement of $\mathcal{O}$. If $\mathcal{D}$ is minimal
in the class of decompositions by connected manifolds, then it cannot be a
refinement of a different decomposition. Hence, $\mathcal{D}=\mathcal{O}$.
\end{proof}

\section{The Whitney Conditions}

In his analysis of stratifications, Whitney introduced two conditions on a
triple of $(M,M^{\prime },x)$ of $C^{1}$-submanifolds of a manifold $W$, and 
$x\in M^{\prime }$, \cite{whitney}. Stratifications satisfying Whitney's
conditions A are called Whitney A stratifications. Similarly,
stratifications satisfying Whitney's conditions B are called Whitney B
stratifications. 

Since we are dealing here with subcartesian spaces, we assume that $M$ and $%
M^{\prime }$ are $C^{\infty }$-submanifolds of $\mathbb{R}^{N}$, and $%
M^{\prime }$ is in the closure $\bar{M}$ of $M$. Let $(x_{n})$ be sequence
of points in $M\in \mathcal{D}$ converging to $x\in M^{\prime }\in \mathcal{D%
}$ such that the sequence of tangent spaces $T_{x_{n}}M$ converges to a
space $D$ in the Grassmannian of $m$-planes in $R^{N}$, where $m=(\dim M).$

\begin{condition}[Whitney Condition A]
$T_{x}M^{\prime }\subseteq D.$
\end{condition}

\begin{condition}[Whitney Condition B]
If $y_{n}$ is a sequence of points in $M^{\prime }$ converging to $x$ and
the sequence of lines $\left\langle x_{n}.y_{n}\right\rangle $ converges to
a line $L$ through $x$, then $L\subset D$.
\end{condition}

Let $\boldsymbol{e}=(e_{1},...,e_{N})$ be the canonical basis of $\mathbb{R}%
^{N}$. Each orthonormal basis $\boldsymbol{b}=(b_{1},...,b_{N})$ in $R^{N}$
is of the form $\boldsymbol{b}=\boldsymbol{e}A$, for a unique $A\in O(N)$.
An orthonormal basis $\boldsymbol{b}=(b_{1},...,b_{N})$ in $\mathbb{R}^{N}$
is said to be adapted to an $m$-dimensional subspace $D$ if the first $m$
vectors $(b_{1},...,b_{m})$ in $\boldsymbol{b}$ form a basis of $D$. The
class of all bases of $\mathbb{R}^{N}$ adapted to $D$ is given by an element 
$\gamma \in O(N)/(O(N)\times O(N-m))$. Using the bijection $A\mapsto 
\boldsymbol{b}=\boldsymbol{e}A$ between $O(N)$ and the space of orthonormal
bases on $\mathbb{R}^{N}$, one can identify the set of all $m$-dimenional
subspaces of $R^{N}$ with the Grassmannian $O(N)/O(N)\times O(N-m)$.

Let $D_{n}$ be a sequence of $m$-dimensional subspaces of $\mathbb{R}^{N}$.
For each $n$, we denote by $\gamma _{n}\in O(N)/(O(N)\times O(N-m))$ the
class of bases in $R^{N}$ adapted to $D_{n}.$ The seqeuence of subspaces $%
D_{n}$ is said to be convergent to an $m$-dimensional subspace $D$ if the
sequence $\gamma _{n}$ converges in $O(N)/(O(N)\times O(N-m))$ to $\gamma $
representing the class of all bases in $D$.

Assume that $D_{n}$ converges to $D$. For each $n$, we can choose a matrix $%
A_{n}\in O(N)$ such that $\boldsymbol{b}_{n}=\boldsymbol{e}A_{n}$ is adapted
to $D_{n}$. Since $O(N)$ is compact, there exists a convergent subsequence $%
A_{n_{k}}$. Let $A=\lim_{k\rightarrow \infty }A_{n_{k}}$. Then $\boldsymbol{b%
}=\boldsymbol{e}A=\lim_{k\rightarrow \infty }\boldsymbol{b}_{n_{k}}$. If $%
b=(b_{1},...b_{N})$ then, for each $i=1,...,N$, the sequence $b_{n_{k}i}$ of 
$i$'th vectors in $\boldsymbol{b}_{n_{k}}=(b_{n_{k}1},...,b_{n_{k}N})$
converges to $b_{i}$.

Let $u_{n}\in D_{n}$ be a convergent sequence of vectors in $\mathbb{R}^{N}$
and $u=\lim_{n\rightarrow \infty }u_{n}$. For each $n$, we can express $%
u_{n} $ in terms of the basis $\boldsymbol{b}_{n}=(b_{n1},...,b_{nN})$
obtaining $u_{n}=a_{n}^{1}b_{n1}+...+a_{n}^{N}b_{nN}$. Since the basis $%
\boldsymbol{b}_{n}$ is orthonormal, for each $i=1,...,N$, we have $%
a_{n}^{i}=u_{n}\cdot b_{ni}$, where $\cdot $ denotes the canonical scalar
product in $\mathbb{R}^{N}$. Hence, 
\[
\lim_{k\rightarrow \infty }a_{n_{k}}^{i}=\lim_{k\rightarrow \infty
}(u_{n_{k}}\cdot b_{n_{k}i})=(\lim_{k\rightarrow \infty }u_{n_{k}})\cdot
(\lim_{k\rightarrow \infty }b_{n_{k}i})=u\cdot b_{i}. 
\]%
Therefore, $u=(u\cdot b_{1})b_{1}+...+(u\cdot b_{N})b_{N}\in D$ because $%
\boldsymbol{b}=(b_{1},...,b_{N})$ is a basis in $D$.

Conversely, if $u=a^{1}b_{1}+....+a^{N}b_{N}$ is a vector in $D$, then $%
u_{n}=a^{1}b_{n_{k}1}+...+a^{N}b_{n_{k}N}$ is in $D_{n_{k}}$ because $%
\boldsymbol{b}_{n_{k}}=(b_{n_{k}1},...,b_{n_{k}N})$ is a basis in in $%
D_{n_{k}}$. Moreover, the sequence $u_{n_{k}}\in D_{n_{k}}$ converges to $u$%
. Hence we have justified the following observation:

\begin{remark}
\label{remark1}Suppose that a sequence $D_{n}$ of $m$-dimensional subspaces
of $R^{N}$ converges to an $m$-dimensional subspace $D$. Then, every
convergent sequence of vectors $u_{n}\in $ $D_{n}$ has a limit in $D$ and
vectors in $D$ are limits of convergent sequences of vectors in $D_{n}$.
\end{remark}

\begin{proposition}
The partition of a subcartesian space $P$ by the family $\mathcal{O}$ of
orbits of all vector fields on $P$ satisfies Whitney's condition\ A.
\end{proposition}

\begin{proof}
Let $O^{\prime }$ and $O$ be orbits in $\mathcal{O}$, $x\in O^{\prime }\cap 
\bar{O}$, and $(x_{n})$ be sequence of points in $M$ converging to $x$ such
that the sequence of tangent spaces $T_{x_{n}}O$ converges to $D$ in the
Grassmannian of $m$-planes in $R^{N}$, where $m=(\dim M).$ First, we need to
show that, if a sequence $(x_{n})$ of points in $O$ converges to $x\in
O^{\prime }$ such that the spaces $T_{x_{n}}O$ converge to $D_{x}\subseteq
T_{x}P$, then $T_{x}O^{\prime }\subseteq D_{x}$. Since $P$ is subcartesian,
we may assume without loss of generality that $x$ has a neighbourhood $U$ in 
$P$ that can be identified with a subset of $\mathbb{R}^{N}$. Each $%
T_{x_{n}}O$ can be identified with the corresponding $m$-dimensional
subspace $D_{n}$ of $\mathbb{R}^{N}$, where $m=\dim O$. Similarly, we
identify $D_{x}\subseteq T_{x}P$ with an $m$-dimensional subspace $D$ of $%
\mathbb{R}^{N}$. By assumption, the sequence $D_{n}$ converges to $D$.

Let $k=\dim O^{\prime }$, $m=\dim O$, and let $%
X_{1},...,X_{k},X_{k+1,}...,X_{m}$ be a family of vector fields on $P$ such
that $X_{1}(x),...,X_{k}(x)$ is a basis for $T_{x}O^{\prime }$ and $%
X_{1},...,X_{k},X_{k+1,}...,X_{m}$ give rise to a frame in $T(O\cap U)$, for
some neighbourhood $U$ of $x$ in $P$. Without loss of generality, we may
assume that $U$ is the neighbourhood of the preceding paragraph and that all
points of the sequence $x_{n}$ in $O$ converging to $x$ are contained in $%
O\cap U$.

For each $i=1,...,m$, the vector field $X_{i}$ is continuous so that $%
X_{i}(x)=\lim_{n\rightarrow \infty }X_{i}(x_{n})$. Since $%
(X_{1}(x),...,X_{k}(x))$ is a frame for $T_{x}O^{\prime }$, it follows that
every vector $u\in T_{x}O^{\prime }$ is of the form $%
u=a^{1}X_{1}(x)+...+a^{k}X_{k}(x)$. Let $%
u_{n}=a^{1}X_{1}(x_{n})+...+a^{k}X_{k}(x_{n})\in T_{x_{n}}O$. Then $%
u=\lim_{n\rightarrow \infty }u_{n}$, and Remark \ref{remark1} implies that $%
u\in D_{x}$. Hence, $T_{x}O^{\prime }\subseteq D_{x}$, which implies
Whitney's condition $A$.
\end{proof}

In general, the family $\mathcal{O}$ of orbits of all vector fields on a
subcartesian space $P$ need not satisfy Whitney's Condition B.

\begin{example}[Spiral]
Let $S$ be the closure of the spiral defined by $r=e^{-\theta }$ in $\mathbb{%
R}^{2}.$ That is $S=S_{0}\cup S_{1},$where $S_{0}=\{(0,0)\}$and $%
S_{1}=\{(e^{-\theta }\cos \theta ,e^{-\theta }\sin \theta )\mid \theta \in 
\mathbb{R}\}.$The slope of $S_{1}$ at $\theta $ is 
\[
m_{\theta }=\frac{(e^{-\theta }\sin \theta )^{\prime }}{(e^{-\theta }\cos
\theta )^{\prime }}=\frac{-e^{-\theta }\sin \theta +e^{-\theta }\cos \theta 
}{-e^{-\theta }\cos \theta -e^{-\theta }\sin \theta }.
\]%
The sequence of points $\mathbf{x}_{n}=(e^{-2\pi n}\cos (2\pi n),e^{-2\pi
n}(\sin 2\pi n))=(e^{-2\pi n},0)$ converges to the origin. Moreover, the
slope of $S_{1}$ at $\mathbf{x}_{n}$ is $m_{2\pi n}=-1$, which implies that
the sequence $T_{\mathbf{x}_{n}}S_{1}$ converges to a line $y=-x.$

For each $n$, the line $L_{n}$ joining $\boldsymbol{x}_{n}$ to the origin $%
\boldsymbol{0}=(0,0)$ is the $y$-axis. Hence, the sequence $L_{n}$ converges
to the $y$-axi\.{s} which is not contained in $\lim_{n\rightarrow \infty }T_{%
\mathbf{x}_{n}}S_{1}$. Thus, our spiral does not satisfies Whitney's
condition $B$.
\end{example}

Nevertheless, there are several Whitney B stratifications which are given by
the family $\mathcal{O}$ of orbits of vector fields on a subcartesian space.

\begin{example}[Whitney's cusp]
Whitney's cusp $S\subseteq \mathbb{R}^{3}$ is the zero level set of $%
F(x,y,z)=y^{2}+x^{3}-z^{2}x^{2}.$ In other words, 
\[
S=\{(x,y,z)\in \mathbb{R}^{3}\mid y^{2}+x^{3}-z^{2}x^{2}=0\}.
\]%
Since $F\in C^{\infty }(\mathbb{R}^{3})$, Implicit function Theorem implies
that $S$ is a smooth manifold in a neighbourhood of every point $(x,y,z)$ in 
$S$ such that $DF(x,y,z)\neq 0$. But 
\[
DF(x,y,z)=(3x^{2}-2xz^{2})dx+2ydy-2zx^{2}dz.
\]%
Hence, $DF(x,y,z)=0$ on the $z$-axis 
\[
S_{1}=\{(x,y,z)\in \mathbb{R}^{3}\mid x=y=0\},
\]%
and $S_{2}=S\backslash S_{1}$ is a smooth manifold. The Hessian of $F$ is 
\[
D^{2}F(x,y,z)=(6x-z^{2})dx^{2}+2dy^{2}-2x^{2}dz-4xzdxdz.
\]%
It has rank $2$ on%
\[
S_{1}^{\pm }=\{(x,y,z)\in \mathbb{R}^{3}\mid x=y=0,\text{ }\pm z>0\},
\]%
and rank $1$ at the origin $S_{0}=\{(0,0,0)\}$.

The function $F$ is invariant under the action of $\mathbb{R}$ on $\mathbb{R}%
^{3}$ given by 
\[
\Phi :\mathbb{R\times R\rightarrow R}:(x,y,z)\mapsto (e^{2t}x,e^{3t}y,e^{t}z)
\]%
generated by a vector field 
\[
X=2x\frac{\partial }{\partial x}+3y\frac{\partial }{\partial y}+z\frac{%
\partial }{\partial z}.
\]%
This action is transitive on $S_{1}^{+}$ and $S_{1}^{-}$. Moreover, since
the Hessian of a smooth function on a manifold is well defined on the set of
critical points of the function, it follows that every diffeomorphism of $%
\mathbb{R}^{3}$ to itself which leaves $F$ invariant preserves the origin.
This implies that the decomposition $S=S_{2}\cup S_{1}^{+}\cup S_{2}^{-}\cup
S_{0}$ is the partition of given by the family of orbits of all smooth
vector fields on $S$. It is of interest to note that this partition is a
stratification of $S$ satisfying Whitney's conditions A and B.
\end{example}

\section{Orbits of a proper group action on a manifold.}

In this section we prove that the stratification by orbit type of the space
of orbits of a proper action of a Lie group $G$ on a manifold $P$ is given
by the family $\mathcal{O}$ of orbits of all vector fields on the orbit
space $R=P/G$ with the differential structure $C^{\infty }(R)$ given by the
ring $C^{\infty }(P)^{G}$ of $G$-invariant smooth functions on $P$. This
results shows that the ring $C^{\infty }(P)^{G}$ encodes information about
the stratification structure of the orbit space $P/G$.

We consider here a proper action 
\[
\Phi :G\times P\rightarrow P:(g,p)\mapsto \Phi (g,p)\equiv \Phi
_{g}(p)\equiv gp 
\]%
of a connected Lie group on a manifold $P$. Properness of $\Phi $ means
that, for every convergent sequence $(p_{n})$ in $P$ and a sequence $(g_{n})$
in $G$ such that the sequence $(g_{n}p_{n})$ is convergent, the sequence $%
(g_{n})$ has a convergent subsequence $(g_{n_{k}})$ and 
\[
\lim_{k\rightarrow \infty }(g_{n_{k}}p_{n_{k}})=\left( \lim_{k\rightarrow
\infty }g_{n_{k}}\right) \left( \lim_{k\rightarrow \infty }p_{n_{k}}\right)
. 
\]

For $p\in P$, the orbit of $G$ through $p$ is the set 
\[
Gp=\{gp\mid g\in G\}.
\]%
Let $R=P/G$ denote the space of $G$-orbits in $P$ with the quotient topology
and let $\rho :P\rightarrow R:p\mapsto Gp$ be the canonical projection.
Since the action $\Phi $ is proper, the orbit space $R$ is a subcartesian
space with the ring $C^{\infty }(R)$ of smooth functions on $R$ given by%
\[
C^{\infty }(R)=\{f\in C^{0}(R)\mid \rho ^{\ast }f\in C^{\infty }(P)\},
\]%
and the projection map $\rho :P\rightarrow R$ is smooth \cite%
{cushman-sniatycki}.

The orbit space $R=P/G$ of a proper action of a Lie group is stratified by
orbit type . Since $\Phi $ is proper, for each $p\in P$,\ the isotropy group 
\[
G_{p}=\{g\in G\mid gp=p\} 
\]%
of $p$ is compact. For each compact subgroup $H\subseteq G$, 
\[
P_{H}=\{p\in P\mid G_{p}=H\} 
\]%
is the set of all points in $P$ of isotropy type $H$. Similarly, 
\[
P_{(H)}=\{p\in P\mid G_{p}\text{ is congruent to }H\mathbf{\}} 
\]%
is the set of all points in $P$ of orbit type $H$. Both $P_{H}$ and $P_{(H)}$
are local submanifolds of $P.$ This means that connected components of $%
P_{H} $ and $P_{(H)}$ are submanifolds of $P$. Connected components of the
projection of $P_{(H)}$ to $R$ are smooth manifolds. They are strata of the
stratification of $R$ by orbit type. For more details see \cite%
{duistermaat-kolk}.

\begin{theorem}
\label{theorem4}The stratification of $R=P/G$ by orbit type coincides with
the partition of $R$ by the family $\mathcal{O}$ of orbits of all vector
fields on \thinspace $R$.
\end{theorem}

\begin{proof}
Theorem \ref{theorem3} implies that it suffices to prove that the
stratification of $R$ by orbit types is minimal and it admits local
extensions of vector fields, Minimality of $R$ has been proved by Bierstone 
\cite{bierstone75}, \cite{bierstone80}. See also Duistermaat \cite%
{duistermaat}. Hence, it remains to prove that the orbit type stratification
of $R$ admits local extensions of vector fields.

Let $M$ be a stratum of the stratification of $R$ by orbit type and $X_{M}$
a smooth vector field on $M$. We want to show that, for each $x\in M$, there
exists a neighbourhood $V\subseteq M$ and a vector field $X$ on $R$ such
that the restrictions to $V$ of $X$ and $X_{M}$ coincide.

Since the action of $G$ on $P$ is proper, for each $p\in \rho ^{-1}(x),$
there exists a slice $\Sigma $ for the action of $G$ at $p$. That is, $%
\Sigma $ is a submanifold of $P$ containing $p$, invariant under the action
of $G_{p},$ and satisfying the following conditions 
\begin{eqnarray}
T_{p}P &=&T_{p}\Sigma \oplus T_{p}(Gp),  \label{1} \\
T_{p^{\prime }}P &=&T_{p^{\prime }}\Sigma +T_{p^{\prime }}(Gp^{\prime })%
\text{ for all }p^{\prime }\in \Sigma ,  \label{2} \\
\text{For }p^{\prime } &\in &\Sigma \text{ and }g\in G\text{, if }gp^{\prime
}\in \Sigma \text{, then }g\in G_{p}.  \label{3}
\end{eqnarray}%
Given a slice $\Sigma ,$ the set 
\[
\tilde{W}=\bigcup_{p^{\prime }\in \Sigma }Gp^{\prime }=\{gp^{\prime
}\mid p^{\prime }\in \Sigma \text{, }g\in G\} 
\]%
is a $G$-invariant neighbourhood $\tilde{W}$ of $Gp$ in $P$. Its projection $%
W=\rho (\tilde{W})$ to $R$ is an open neighbourhood of $x$ in $R$.

Let $\tilde{M}$ be the connected component of $P_{G_{p}}$ that contains $p$.
As we have said before, $\tilde{M}$ is a submanifold of $P$. Moreover, $%
M=\rho (\tilde{M})$. The intersection $\tilde{M}\cap W$ is an open
submanifold of $\tilde{M}$.

\textbf{Lemma.}  $\tilde{M}\cap \Sigma $ is a submanifold of $\Sigma $
diffeomorphic to $M\cap W$.

The condition (\ref{3}) states that, if $p^{\prime }$ and $gp^{\prime }$ are
in $\Sigma $, then $g\in G_{p}$. Moreover, $p^{\prime }$ and $gp^{\prime }$
in $\tilde{M}$, implies that $G_{p^{\prime }}=G_{gp^{\prime }}=G_{p}$.
Hence, $g\in G_{p}$ implies that $g\in G_{p^{\prime }}$ and $gp^{\prime
}=p^{\prime }$. Thus, $\tilde{M}\cap \Sigma $ intersects fibres of the
projection map $\rho :P\rightarrow R$ in at most single points. Therefore,
the restriction $\mu $ of the projection map $\rho $ to $\tilde{M}\cap
\Sigma $ is a bijection of $\tilde{M}\cap \Sigma $ onto $M\cap W$. To show
that $\mu :\tilde{M}\cap \Sigma \rightarrow M\cap W$ is a diffeomorphism it
suffices to show that $\mu $ and $\mu ^{-1}$ are smooth.

The space $C^{\infty }(\tilde{M}\cap \Sigma )$ is generated by restrictions
to $\tilde{M}\cap \Sigma $ of smooth functions on $P$. On the other hand, $%
C^{\infty }(M\cap W)$ is generated by restrictions to $M\cap W$ of functions
in $C^{\infty }(R)=\{\rho _{\ast }f\mid f\in C^{\infty }(P)^{G}\}$.

First, we show that $\mu :\tilde{M}\cap \Sigma \rightarrow M\cap W$ is
smooth. Consider a function $f_{M\cap W}\in C^{\infty }(M\cap W)$. We need
to show that $\mu ^{\ast }f_{M\cap W}\in C^{\infty }(\tilde{M}\cap \Sigma ).$
For each point $x^{\prime }\in M\cap W$, there exists a neighbourhood $U$ of 
$x^{\prime }$ in $M\cap W$ and a function $\rho _{\ast }f\in C^{\infty }(R)$
such that the restriction of $\rho _{\ast }f$ to $U$ coincides with the
restriction to $U$ of $f_{M\cap W}$. For each $p^{\prime \prime }\in \mu
^{-1}(U)$, we have 
\[
\mu ^{\ast }f_{M\cap W}(p^{\prime \prime })=f_{M\cap W}(\mu (p^{\prime
\prime }))=\rho _{\ast }f(\mu (p^{\prime \prime }))=f(p^{\prime \prime
})=f_{\mid \tilde{M}\cap \Sigma }(p^{\prime \prime }),
\]%
since $f$ is $G$-invariant. Hence, $\mu ^{\ast }f_{M\cap W}$ restricted to $%
\mu ^{-1}(U)$ coincides with the restriction to $\mu ^{-1}(U)$ of $f\in
C^{\infty }(P)^{G}$. Since this result is valid for each $x^{\prime }\in
M\cap W,$ it follows that $\mu ^{\ast }f_{M\cap W}\in C^{\infty }(\tilde{M}%
\cap \Sigma ).$ However, $f_{M\cap W}$ is an arbitrary smooth function on $%
M\cap W$. Hence, $\mu :\tilde{M}\cap \Sigma \rightarrow M\cap W$ is smooth.

Next, we want to show that $\mu ^{-1}:M\cap \Sigma \rightarrow \tilde{M}\cap
W$ is smooth. Consider a function $f_{\tilde{M}\cap \Sigma }$ in $C^{\infty
}(\tilde{M}\cap \Sigma )$. We need to show that $f_{\tilde{M}\cap \Sigma
}\circ \mu ^{-1}$ is in $C^{\infty }(M\cap W)$. Given $q\in \tilde{M}\cap
\Sigma ,$ there exists a compactly supported function $f_{\Sigma }$ on $%
\Sigma $ which coincides with $f_{\tilde{M}\cap \Sigma }$ on a neighbourhood 
$\tilde{U}$ of $q$ in $\tilde{M}\cap \Sigma $. Let $\bar{f}_{\Sigma }$ be
the $G_{p}$ invariant function on $\Sigma $ obtained by averaging $f_{\Sigma
}$ over $G_{p}$. Then, $f_{\Sigma \mid \tilde{U}}=\bar{f}_{\Sigma \mid 
\tilde{U}}$ because $\tilde{U}\subseteq \tilde{M}\cap \Sigma $ and all
points in $\tilde{M}\cap \Sigma $ have isotropy group $G_{p}$. Let $f_{%
\tilde{W}}$ be a function on $\tilde{W}$ defined by%
\begin{equation}
f_{\tilde{W}}(gp^{\prime })=\bar{f}_{\Sigma }(p^{\prime })\text{ for all }%
p^{\prime }\in \Sigma \text{ and }g\in G\text{.}  \label{fW}
\end{equation}%
The function $f_{\tilde{W}}$ is well defined by equation (\ref{fW}) and is $G
$-invariant because $\bar{f}_{\Sigma }$ is $G_{p}$-invariant. If $g^{\prime
}p^{\prime }=g^{\prime \prime }p^{\prime \prime }$, with $p^{\prime }$ and $%
p^{\prime \prime }$ in $\Sigma $, then $p^{\prime }=(g^{\prime
})^{-1}g^{\prime \prime }p^{\prime \prime }$ implies that $(g^{\prime
})^{-1}g^{\prime \prime }\in G_{p}$ and 
\[
f_{\tilde{W}}(g^{\prime }p^{\prime })=\bar{f}_{\Sigma }(p^{\prime })=\bar{f}%
_{\Sigma }((g^{\prime })^{-1}g^{\prime \prime }p^{\prime \prime })=\bar{f}%
_{\Sigma }(p^{\prime \prime })=f_{\tilde{W}}(p^{\prime \prime })\text{.}
\]%
Moreover, for every $g,g^{\prime }\in G$, and $p^{\prime }\in \Sigma ,$ 
\[
f_{\tilde{W}}(g(g^{\prime }p^{\prime }))=f_{\tilde{W}}(gg^{\prime }p^{\prime
})=\bar{f}_{\Sigma }(p^{\prime })=f_{\tilde{W}}(g^{\prime }p^{\prime }).
\]%
Since $\bar{f}_{\Sigma }$ is compactly supported, it follows that there is a
a $G$-invariant open set $\tilde{V}$ in $P$ containing the support of $f_{%
\tilde{W}}$ and such that the closure of $\tilde{V}$ is in $\tilde{W}$.
Therefore, there exists an extension of $f_{\tilde{W}}$ to a smooth function 
$f$ on $P$ that $f$ vanishes on the complement of $\tilde{V}$. Moreover, $%
f\in C^{\infty }(P)^{G}$ because $f_{\tilde{W}}$ is $G$-invariant. For each $%
p^{\prime }\in \tilde{U}$, we have 
\[
f_{\tilde{M}\cap \Sigma }\circ \mu ^{-1}(\rho (p^{\prime }))=f_{\tilde{M}%
\cap \Sigma }(p^{\prime })=f(p^{\prime })=\rho _{\ast }f(\rho (p^{\prime })).
\]%
Hence, $f_{\tilde{M}\cap \Sigma }\circ \mu ^{-1}$ restricted to an open
neighbourhood $U=\rho (\tilde{U})=\mu (\tilde{U})$ of $\rho (q)$ in $M\cap W$
coincides with the restriction to $U$ of $\rho _{\ast }f\in C^{\infty }(R)$.
Since it holds for every point $\rho (q)\in M\cap \Sigma $, it follows that $%
(\mu ^{-1})^{\ast }f_{\tilde{M}\cap \Sigma }=f_{\tilde{M}\cap \Sigma }\circ
\mu ^{-1}\in C^{\infty }(M\cap \Sigma )$. Hence, $\mu ^{-1}:M\cap \Sigma
\rightarrow \tilde{M}\cap \Sigma $ is smooth. This completes the proof of
our lemma.

We continue with the proof of Theorem \ref{theorem4}. Consider a smooth
vector field $X_{M}$ on $M$. For each $x\in M$, consider $p\in \rho ^{-1}(x)$
and let $\tilde{M}$ be the connected component of $\rho ^{-1}(M)$ containing 
$p$. Let $\Sigma $ be a slice at $p$ for the action of $G$ on $P$, and $%
W=\rho (\Sigma )$. We have shown that $\mu :\tilde{M}\cap \Sigma \rightarrow
M\cap W$ is a diffeomorphism Let $f_{M\cap W}$ be a compactly supported
smooth function on $M\cap W$ such that $f_{M\cap W}(x^{\prime })=1$ for all $%
x^{\prime }$ in a neighbourhood of $x$ in $M\cap W$. The product $f_{M\cap
W}X_{M}$ is a vector field on $M\cap W$ which can be pushed forward by $\mu
^{-1}:M\cap W\rightarrow \tilde{M}\cap \Sigma $ to a vector field $\mu
_{\ast }^{-1}(f_{M\cap W}X_{M})$ on $\tilde{M}\cap \Sigma $. Let $c:t\mapsto
c(t)$ be an integral curve of $\mu _{\ast }^{-1}(f_{M\cap W}X_{M})$. Since $%
c(t)$ is contained in $M,$ for each $g\in G_{p}$ we have $gc(t)=c(t)$.
Hence, $c$ is invariant under the action of $G_{p}$. This implies that $\mu
_{\ast }^{-1}(f_{M\cap W}X_{M})$ is $G_{p}$-invariant. That is, for each $%
g\in G_{p}$, 
\[
T\Phi _{g}\circ \mu _{\ast }^{-1}(f_{M\cap W}X_{M})\circ \Phi _{g^{-1}}=\mu
_{\ast }^{-1}(f_{M\cap W}X_{M}).
\]

Since $\mu _{\ast }^{-1}(f_{M\cap W}X_{M})$ is compactly supported in a
neighbourhood of $p=\mu ^{-1}(x)$ in $\tilde{M}\cap \Sigma $, it can be
extended by zero to a vector field $X_{\Sigma }$ on $\Sigma $. Note that $%
X_{\Sigma }$ is $G_{p}$-invariant, since, for each $g\in G_{p}$ and $%
p^{\prime }\in \Sigma ,$ either (1) $p^{\prime }\in M\cap \Sigma $ or (2) $%
p^{\prime }\notin M\cap \Sigma $. If (1) $p^{\prime }\in M\cap \Sigma ,$
then 
\[
T\Phi _{g}\circ X_{\Sigma }\circ \Phi _{g-1}(p^{\prime })=T\Phi _{g}\circ
\mu _{\ast }^{-1}(f_{M\cap W}X_{M})\circ \Phi _{g-1}(p^{\prime })=\mu _{\ast
}^{-1}(f_{M\cap W}X_{M})(p^{\prime })=X_{\Sigma }(p^{\prime }).
\]%
If $(2)$ $p^{\prime }\notin M\cap \Sigma $ and $g^{-1}p^{\prime }\in M\cap
\Sigma $, then $G_{g^{-1}p^{\prime }}=G_{p}$ and $g\in G_{p}$ implies that $%
g\in G_{p}=G_{g^{-1}p^{\prime }}$ so that $g^{-1}p^{\prime
}=g(g^{-1}p^{\prime })=p^{\prime }$ and we have contradiction with the
assumption that $p^{\prime }\notin M\cap \Sigma $. Hence, $g^{-1}p^{\prime
}\notin M\cap \Sigma .$ In this case $X_{\Sigma }(p^{\prime })=X_{\Sigma
}(g^{-1}p^{\prime })=0,$ and%
\[
T\Phi _{g}\circ X_{\Sigma }\circ \Phi _{g-1}(p^{\prime })=T\Phi _{g}\circ
X_{\Sigma }(g^{-1}p^{\prime })=T\Phi _{g}(0)=0=X_{\Sigma }(p^{\prime }).
\]%
In any case, 
\[
T\Phi _{g}\circ X_{\Sigma }(p^{\prime })=X_{\Sigma }(gp^{\prime })
\]%
for all $g\in G_{p}$ and $p^{\prime }\in \Sigma $

We can extend $X_{\Sigma }$ to a $G$-invariant vector field $X_{\tilde{W}}$
on $\tilde{W}$ by setting%
\[
X_{\tilde{W}}(gp^{\prime })=T\Phi _{g}(X_{\Sigma }(p^{\prime })) 
\]%
for every $g\in G$ and $p^{\prime }\in \Sigma $. It is well defined since,
if $gp^{\prime }=g^{\prime }p^{\prime \prime }$ for $p^{\prime },p^{\prime
\prime }\in \Sigma $, then $g^{-1}g^{\prime }\in G_{p}$ and 
\[
T\Phi _{g^{\prime }}(X_{\Sigma }(p^{\prime \prime }))=T\Phi _{g}(T\Phi
_{g^{-1}g^{\prime }}(X_{\Sigma }(p^{\prime \prime }))=T\Phi _{g}(X_{\Sigma
}(\Phi _{g^{-1}g^{\prime }}(p^{\prime \prime })). 
\]%
If $p^{\prime \prime }\in \Sigma \cap \tilde{M}$, then $G_{p^{\prime \prime
}}=G_{p}$ and $\Phi _{g^{-1}g^{\prime }}(p^{\prime \prime })=g^{-1}g^{\prime
}p^{\prime \prime }=p^{\prime \prime }$ and $X_{\tilde{W}}(gp^{\prime })=X_{%
\tilde{W}}(g^{\prime }p^{\prime \prime }).$ If $p^{\prime \prime }\notin
\Sigma \cap \tilde{M}$ and $g^{-1}g^{\prime }\in G_{p}$, then $p^{\prime
}\notin \Sigma \cap \tilde{M}$, and $X_{\tilde{W}}(gp^{\prime })=0=X_{\tilde{%
W}}(g^{\prime }p^{\prime \prime }).$

Finally, we can extend $X_{\tilde{W}}$ to a $G$-invariant vector field $X$
on $P$, by setting $X(p^{\prime })=X_{\tilde{W}}(p^{\prime })$ for $%
p^{\prime }\in \tilde{W}$ and $X(p^{\prime \prime })=0$ for $p^{\prime
\prime }\notin \tilde{W}.$ Since $X$ is $G$-invariant, it restricts to a
derivation of $C^{\infty }(P)^{G}$ which is equivalent to a derivation of $%
C^{\infty }(R)$. This derivation is a vector field on $R$ because it comes
from a vector field on $P$.

Thus, the stratification of $R$ by orbit type admits local extensions of
vector fields. Since it is also minimal, Theorem \ref{theorem3} implies that
it coincides with the partition of $R$ by the family $\mathcal{O}$ of orbits
of all vector fields on $R$.
\end{proof}

\end{document}